# A Note on the Koethe Conjecture

S. K. PANDEY

Department of Mathematics,

SPUP, Vigyan Nagar, Jodhpur, India.

Email: skpandey12@gmail.com

**Abstract**

In this note we provide some results related to the Koethe conjecture and exhibit that the condition '$R$ satisfies the Koethe conjecture' given in [1-2, theorem 2.6] is superfluous at least under certain conditions described in this note.

## 1 Introduction

The well known Koethe conjecture was introduced in 1930. It is a longstanding problem in ring theory, with a long and complicated history. It can be stated as follows. Every one sided nil ideal of a ring $R$ is contained in a two sided nil ideal of $R$ [3].

Let $N(R)$ is the set of all nilpotent elements of a ring $R$. As per theorem 2.6 [1-2], the following statements are equivalent for a ring $R$.

  (a) $N(R)$ is closed under addition.
  (b) $N(R)$ is closed under multiplication and $R$ satisfies Koethe conjecture.
  (c) $N(R)$ is a subring of $R$.

Thus it has been shown in [1-2] that $R$ satisfies the Koethe conjecture provided $N(R)$ is a subring of $R$. Also, it has been shown in [1-2] that if $N(R)$ is multiplicatively closed and $R$ satisfies the Koethe conjecture then $N(R)$ is a subring of $R$.

In this note we exhibit that the condition '$R$ satisfies the Koethe conjecture' given above (theorem 2.6, [1-2]) is superfluous at least under certain conditions described here.

In addition we provide some results when $N(R)$ does not form a subring of $R$. Also, we partially address an open question appeared in [1-2] which asks that if $N(R)$ is multiplicatively closed, then does it imply that $R$ satisfies the Koethe conjecture?

**2 Some Results**

We provide the following results.

**Proposition 1.** Let $R$ is a noncommutative ring and $N(R)$ is the set of all nilpotent elements of $R$. If $R$ fails to satisfy the Koethe conjecture, then $N(R)$ is not a subring of $R$.

**Proof.** Let $R$ is a noncommutative ring and $N(R)$ is the set of all nilpotent elements of $R$. If $N(R)$ is a subring of $R$, then $R$ satisfies the Koethe conjecture [1-2, theorem 2.6]. Therefore it follows that if $R$ fails to satisfy the Koethe conjecture, then $N(R)$ is not a subring of $R$.

**Proposition 2.** The converse of the proposition 1 is not true.

**Proof.** Let $R$ is the ring of all $2 \times 2$ matrices over the field of order two. Then $N(R)$ is not a subring of $R$. However $R$ satisfies the Koethe conjecture.

**Proposition 3.** Let $R$ is a noncommutative ring and $N(R)$ is the set of all nilpotent elements of $R$. Then the following are equivalent.

  i)  $N(R)$ is not a subring of $R$.
  ii) $N(R)$ is not additively closed.

**Proof.** Let $R$ is a noncommutative ring and $N(R)$ is the set of all nilpotent elements of $R$. It is trivial that ii) $\Rightarrow$ i). Let $N(R)$ is additively closed. Then $N(R)$ is a subring of $R$ [1-2, theorem 2.6]. This is a contradiction. Hence if $N(R)$ is not a subring of $R$, then $N(R)$ is not additively closed. Thus i) $\Rightarrow$ ii).

**Proposition 4.** Let $R$ is a noncommutative ring and $N(R)$ is the set of all nilpotent elements of $R$. Then $N(R)$ is not multiplicatively closed is equivalent to $N(R)$ is not additively closed.

**Proof.** Let $N(R)$ is not multiplicatively closed. This implies that $N(R)$ is not a subring of $R$. Now from proposition 3, it follows that $N(R)$ is not additively closed.

**Remark 1.** It may be noted that there is an open question given in [1-2]. The question is as follows. Let $R$ is a ring such that $N(R)$ is multiplicatively closed. Does $R$ satisfy the Koethe conjecture ? It should be emphasized that this open question posed in [1-2] remains open if the converse of the proposition 4 is not true. This is due to the fact that if $N(R)$ is multiplicatively as well as additively closed, then $R$ satisfies the Koethe conjecture.

In the next few results we exhibit that the condition '$R$ satisfies the Koethe conjecture' given in [1-2, theorem 2.6 ] is superfluous at least under following conditions.

**Proposition 5.** Let $R$ is a ring such that $N(R)$ is multiplicatively closed. Then $N(R)$ is additively closed and it is superfluous to assume that $R$ satisfies the Koethe conjecture provided each element of $N(R)$ is of index two.

**Proof.** Let $R$ be a ring such that $N(R)$ is multiplicatively closed. Let $a,b \in N(R)$ are any two elements of $N(R)$. Then it easily follows that $a+b \in N(R)$ since $(a+b)^4 = 0$. Therefore $N(R)$ a subring of $R$ and hence $R$ satisfies the Koethe conjecture [1-2, theorem 2.6]. Clearly, in this case, in order to prove that $N(R)$ is a subring of $R$ it suffices that $N(R)$ is multiplicatively closed and unlike [1-2] one does not require additional condition, namely, $R$ satisfies the Koethe conjecture.

**Theorem 1.** Let $R$ is a ring such that $N(R)$ is multiplicatively closed. Then $N(R)$ is additively closed and it is superfluous to assume that $R$ satisfies the Koethe conjecture provided any one of the following hold.

1) $N(R)[x] \subseteq N(R[x])$.
2) $N(R)[x] = N(R[x])$.

Here $R[x]$ stands for polynomial ring defined over $R$.

**Proof.** Please refer [4, theorem 2.11].

**Theorem 2.** Let $R$ is a ring such that $N(R)$ is multiplicatively closed. Then $N(R)$ is additively closed and it is superfluous to assume that $R$ satisfies the Koethe conjecture provided $N(R)[[x]] = N(R[[x]])$. Here $R[[x]]$ stands for power series ring defined over $R$.

**Proof.** Please refer [4, theorem 2.11].

**Remark 2.** Proposition 4 and the remark following the proof of Proposition 4 lead to the following open question.

**Question.** Is the converse of the Proposition 4 true?